%% file: Z-Maslov-nospringer.tex
\documentclass[11pt]{article}

\usepackage{amsmath}
\usepackage{amsfonts}
\usepackage{amssymb}
\usepackage{amsthm}
\usepackage[english]{babel}
\usepackage[applemac]{inputenc}
\usepackage[all]{xy}
\usepackage{amsthm,amssymb,amsbsy,amsmath,amsfonts,amssymb,amscd,mathrsfs}

\newcommand{\R}{\mathbb{R}}
\newcommand{\C}{\mathbb{C}}
\newcommand{\Z}{\mathbb{Z}}
\newcommand{\RP}{\mathbb{R}\textrm{P}}

\newtheorem{theorem}{Theorem}

\newtheorem{proposition}[theorem]{Proposition}

\theoremstyle{definition} 
\newtheorem*{definition}{Definition}
\theoremstyle{remark}
\newtheorem{remark}{Remark}
\newtheorem{example}[remark]{Example}

\newcommand{\bi}{\begin{itemize}}
\newcommand{\iii}{\item}
\newcommand{\ei}{\end{itemize}}
\newcommand{\bdeff}{\begin{definition}}
\newcommand{\edeff}{\end{definition}}
\newcommand{\bt}{\begin{theorem}}
\newcommand{\et}{\end{theorem}}
\newcommand{\bp}{\begin{proposition}}
\newcommand{\ep}{\end{proposition}}
\newcommand{\bqn}{\begin{equation}}
\newcommand{\eqn}{\end{equation}}
\newcommand{\bex}{\begin{example}}
\newcommand{\eex}{\end{example}}

\newcommand{\tx}{\text}

\newcommand{\U}{U}
\newcommand{\FF}{F}
\newcommand{\JJ}{J}
\newcommand{\Crit}{C}
\newcommand{\Lagr}{\overline{C}}

\newcommand{\cb}{\R^{k}}

\makeindex             

\begin{document}

\title{Geometry of Maslov cycles}
\author{Davide Barilari\footnote{CNRS, CMAP, \'Ecole Polytechnique, Paris and \'Equipe INRIA GECO Saclay-\^Ile-de-France, Paris, France - Email: \texttt{barilari@cmap.polytecnique.fr}} $\,$and Antonio Lerario\footnote{Department of Mathematics, Purdue University - Email: \texttt{alerario@math.purdue.edu}} }
\date{}
%
%
\maketitle

\begin{abstract}We introduce the notion of \emph{induced Maslov cycle}, which describes and unifies geometrical and topological invariants of many apparently unrelated problems, from Real Algebraic Geometry to sub-Riemannian Geometry.\end{abstract}

\section{Introduction}
In this paper, dedicated to Andrei A. Agrachev in the occasion of his 60th birthday, we survey and develop some of his ideas on the theory of quadratic forms and its applications, from real algebraic geometry  to the study of second order conditions in optimal control theory. The investigation of these problems and their geometric interpretation in the language of symplectic geometry is in fact one of the main contribution of Agrachev's research of the 80s-90s (see for instance \cite{Agrachev1, Agrachev2, AgrachevLerario}) and these techniques are still at the core of his more recent research (see the forthcoming preprints \cite{AgrachevQuHo,AgCurv}).

Also, this survey can be  interpreted as an attempt of the authors to give a unified presentation of the two a priori unrelated subjects of their dissertations under Agrachev's supervision, namely sub-Riemannian geometry and the topology of sets defined by quadratic inequalities. The unifying language comes from symplectic geometry and uses the notion of \emph{Maslov cycle}, as we will discuss in a while.

To start with we introduce some notation. The set $L(n)$ of all $n$-dimensional Lagrangian subspaces of $\R^{2n}$ (with the standard symplectic structure) is called the \emph{Lagrangian Grassmannian}; it is a compact submanifold of the ordinary Grassmannian and once we fix one of its points $\Delta$, we can consider  the algebraic set: 
$$\Sigma=\{\Pi\in L(n)\,|\, \Delta \cap \Pi \neq 0\},$$
(this is what is usually referred to as the \emph{train of $\Delta$}, or the \emph{universal Maslov cycle}).\\
The main idea of this paper is to study \emph{generic} maps $f:X\to L(n)$, for $X$ a smooth manifold, and the geometry of the preimage under $f$ of the cycle $\Sigma$. Such a preimage $f^{-1}(\Sigma)$ is what we will call the \emph{induced} Maslov cycle.  

It turns out that many interesting problems can be formulated in this setting and our goal is to describe a kind of \emph{duality} that allows to get geometric information on the map $f$  by replacing its study with the geometry of $f^{-1}(\Sigma)$. 

To give an example, the Maslov cycle already provides information on the topology of $L(n)$ itself. In fact $\Sigma$ is a \emph{cooriented algebraic hypersurface} smooth outside a set of codimension \emph{three} and its intersection number with a generic map $\gamma:S^1\to L(n)$ computes $[\gamma]\in \pi_1(L(n))\simeq \Z$.

The theory of quadratic forms naturally appears when we look at the local geometry of the Lagrangian Grassmannian: in fact $L(n)$ can be seen as a compactification of the space $Q(n)$ of real quadratic forms in $n$ variables and, using this point of view, the Maslov cycle $\Sigma$ is a compactification of the space of \emph{degenerate} forms. 

Given $k$ quadratic forms $q_1, \ldots, q_k$ we can construct the map:
$$f:S^{k-1}\to L(n),\quad (x_1, \ldots, x_k)\mapsto x_1q_1+\cdots+x_k q_k.$$
In fact the image of this map is contained in the \emph{affine} part of $L(n)$ and its homotopy invariants are trivial. Neverthless the induced Maslov cycle $f^{-1}(\Sigma)$ has a nontrivial geometry and can be used to study the topology of:
$$X=\{[x]\in \RP^{n-1}\,|\, q_1(x)=\cdots=q_k(x)=0\}.$$
More specifically, it turns out that as a first approximation for the topology of $X$ we can take the ``number of holes"  of $f^{-1}(\Sigma)$. Refining this approximation procedure amounts to exploit how the coorientation of $\Sigma$ is pulled-back by $f$.

In some sense this is the idea of the study of (locally defined) families of quadratic forms and their degenerate locus, and the set of Lagrange multipliers for a variational problem admits the same description. In fact one can consider two smooth maps between manifolds $F:U\to M$ and $J:U\to \R$ and ask for the study of critical points of $J$ on level sets of $F$. 
With this notation the manifold of \emph{Lagrange multipliers} is defined to be:
$$C_{F,J}=\{(u, \lambda)\in F^*(T^*M)\, |\, \lambda D_uF-d_uJ=0\}.$$

Attached to every point $(u,\lambda)\in C_{F,J}$ there is a quadratic form, namely the Hessian of $J|_{F^{-1}(F(u))}$ evaluated at $u$, and using this family of quadratic forms we can still define an induced Maslov cycle $\Sigma_{F,J}$ (the definition we will give in the sequel is indeed more intrinsic). 

This abstract setting includes for instance the geodesic problem in Riemannian and sub-Riemannian geometry (and even more general variational problems). In this case  the set $U$ parametrizes the space of admissible curves, $F$ is the \emph{end-point map} (i.e. the map that assigns to each admissible curve its final point), and $J$ is the \emph{energy} of the curve. The problem of finding critical points of the energy on a fixed level set of $F$ corresponds precisely to the geodesic problem between two fixed points on the manifold $M$.

In this context $\Sigma_{F,J}$ corresponds to points where the Hessian of the energy is degenerate and its geometry is related to the structure of \emph{conjugate locus} in sub-Riemannain geometry. Moreover the way this family of quadratic forms (the above mentioned Hessians) degenerates translates into optimality properties of the corresponding geodesics.

Rather than a systematic and fully detailed treatment we try to give the main ideas, giving only some sketches of the proofs (providing references where possible) and offering a maybe different perspective in these well-estabilished research fields. 

Our presentation is strongly influenced by the deep insight and the ideas of A. A. Agrachev. We are extremely grateful to him for having shown them, both in mathematics and in life, the elegance of simpleness.

\section{Lagrangian Grassmannian and universal Maslov cycles}
\subsection{The Lagrangian Grassmannian}

Let us consider $\R^{2n}$ with its standard symplectic form $\sigma$. A vector subspace $\Lambda $ of $\R^{2n}$ is called 
\emph{Lagrangian} if it has dimension $n$ and $\sigma|_{\Lambda}\equiv 0.$  The \emph{Lagrange Grassmannian} $L(n)$ in $\R^{2n}$ is the set of its $n$-dimensional Lagrangian subspaces.

\bp $L(n)$ is a  compact submanifold of the Grassmannian of $n$-planes in $\R^{2n}$; its dimension is $n(n+1)/2.$
\ep
Consider indeed the set $\Delta^{\pitchfork}=\{\Lambda \in L(n)\,|\,  \Lambda \cap \Delta =0 \}$ of all lagrangian subspaces that are transversal to a given one $\Delta\in L(n)$.
Clearly $\Delta^{\pitchfork}\subset L(n)$ is an open subset and
\begin{equation}\label{cover}L(n)=\bigcup_{\Delta \in L(n)} \Delta^{\pitchfork}.\end{equation}
It is then sufficient to find some coordinates on these open subsets. Let us fix a Lagrangian complement $\Pi$ of $\Delta$ (which always exists but is not unique).
Every $n$-dimensional subspace $\Lambda \subset \R^{2n}$ which is transversal to $\Delta$ is the graph of a linear map from $\Pi$ to $\Delta$. Choosing a basis on $\Delta$ and $\Pi$, this linear map is represented in coordinates by a matrix $S$ such that: 
$$\Lambda \cap \Delta=0 \Leftrightarrow \Lambda=\{(x,S x), \,x \in \Pi\simeq \mathbb{R}^{n}\}.$$
 Hence the open set $\Delta^{\pitchfork}$ of all Lagrangian subspaces that are transversal to $\Delta$ is parametrized by the set of symmetric matrices, that gives coordinates on this open set.
This also proves that the dimension of $L(n)$ is $n(n+1)/2$. Notice finally that, being $L(n)$ a closed set in a compact manifold, it is itself compact.

Fix now an element $\Lambda \in L(n)$. The tangent space $T_{\Lambda }L(n)$ to the Lagrange Grassmannian at the point $\Lambda$ can be canonically identified with set of quadratic forms on the space $\Lambda$ itself: 
$$T_{\Lambda}L(n)\simeq Q(\Lambda).$$
Indeed consider a smooth curve $\Lambda(t)$ in $L(n)$ such that $\Lambda(0)=\Lambda$ and denote by $\dot{\Lambda}\in T_{\Lambda}L(n)$ its tangent vector. For any point $x\in \Lambda$ and any smooth extension $x(t)\in \Lambda(t)$ we define the quadratic form:
$$\dot{\Lambda}: x \mapsto \sigma(x,\dot{x}),\qquad \dot{x}= \dot{x}(0).$$
An easy computation shows that this is indeed well defined; moreover writing $\Lambda(t)=\{(x,\,S(t)x),\,x\in \R^{n}\}$ then the quadratic form  $\dot{\Lambda}$ associated to the tangent vector of $\Lambda(t)$ at zero is represented by the matrix $\dot S(0)$, i.e.
$\dot{\Lambda}(x)=x^{T}\dot{S}(0)x$.

We stress that this representation using \emph{symmetric} matrices works only for coordinates induced by a \emph{Lagrangian} splitting $\R^{2n}=\Pi\oplus \Delta,$ i.e. $\Pi$ and $\Delta$ are both lagrangian.
\bex[The Lagrange Grassmannians $L(1)$ and $L(2)$]
Since every line in $\R^2$ is Lagrangian (the restriction of a skew-symmetric form to a one-dimensional subspace must be zero), then $L(1)\simeq \RP^1$.\\
The case $n=2$ is more interesting. We first notice that each 2-plane $W$ in $\R^4$ defines a unique (up to a multiple) degenerate 2-form $\omega$ in $\Lambda^2 \R^4$, by $W=\ker \omega$. Thus there is a map:
$$p:G(2,4)\to\mathbb{P}(\Lambda^2 \R^4)\simeq \RP^5.$$
This map is called the \emph{Plücker embedding}; its image is a projective quadric of signature $(3,3).$ The restriction of $p$ to $L(2)$ maps to 
$$p(L(2))=\{[\omega]\,|\,\ker \omega\neq 0 \quad \textrm{and}\quad \omega \wedge \sigma =0\},$$ which is the intersection of the image of $p$ with an hyperplane in $\RP^5$, i.e. the zero locus of the restriction of the above projective quadric to such hyperplane. In particular $L(2)$ is diffeomorphic to a smooth quadric of signature $(2,3)$ in $\RP^4$.
\eex

 \subsection{Topology of Lagrangian Grassmannians}
It is possible to realize the Lagrange Grassmannian as a homogeneous space, through an action of  the unitary group $U(n)$. In fact we have a homomorphism of groups $\phi:GL(n, \C)\to GL(2n, \R)$ defined by:
$$\phi: A+iB\mapsto \left(\begin{matrix} A &B\\ -B &A\end{matrix} \right),$$ and the image of the unitary group is contained in the symplectic one. In particular for every lagrangian subspace $\Lambda \subset \R^{2n}$ and every $M$ in $U(n)$ the vector space $\phi(M)\Lambda$ is still lagrangian. This defines the action of $U(n)$ on $L(n)$; the stabilizer of a point is readily verified to be the group $O(n)$ and we get:
$$L(n)\simeq U(n)/O(n).$$ 
The cohomology of $L(n)$ can be studied applying standard techniques to the fibration $U(n)\to L(n)$ and working with $\Z_2$ coefficients\footnote{Unless differently stated, all homology and cohomology groups will be with  $\Z_2$ coefficients.} we get a ring isomorphism $H^*(L(n))\simeq H^*(S^1\times \cdots \times S^n)$; we refer the reader to \cite{Fuks} for more details.\\
For our purposes we will need an explict description of the fundamental group of $L(n)$ and this can be obtained as follows. We first consider the map $\det^2:U(n)\to S^1$ (the square of the determinant). Multiplication by a matrix of $O(n)$ does not change the value of the square of the determinant, thus we get a surjective map
$${\det}^2 :L(n)\to S^1.$$
This map also is a fibration and with simply connected fibers, each one of them being diffeomorphic to $SU(n)/SO(n)$. In particular it follows that it realizes an isomorphism of fundamental groups and:
$$\pi_1(L(n))\simeq \Z.$$

\subsection{The universal Maslov cycle}
Since the fundamental group of $L(n)$ is $\Z$, then the $1$-form $d\theta/2\pi$ on $S^1$ (the class of this form generates its first cohomology group with integer coefficients) pulls-back via $\det^2$ to a $1$-form on $L(n)$ whose cohomology class $\mu$ generates $H^1(L(n), \Z)$:
$$\mu=\bigg[\frac{1}{2\pi}({\det}^2)^*\ d\theta\bigg] \in H^1(L(n), \Z).$$
Such a class is usually referred to as the \emph{universal Maslov class} (see \cite{ArnoldChCl, ArnoldGivental}). Once we fix a Lagrangian space $\Delta \in L(n)$ it is possible to define a cooriented algebraic cycle in $L(n)$ which is Poincaré dual to $\mu$; such cycle is called the \emph{train of $\Delta$} and is defined as follows:
$$\Sigma_{\Delta}=\{\Lambda \in L(n) \, |\, \Lambda \cap \Delta \neq 0\}=L(n)\backslash \Delta^{\pitchfork}.$$
Here the subscript denotes the dependence on $\Delta$ and when no confusion arises we will omit it: a different choice of $\Delta$ produces an homologous train (in fact just differing by a symplectic transformation). We will discuss the geometry of $\Sigma$ in greater detail in the next section; what we need for now is that this is an \emph{algebraic} hypersurface whose singularities have codimension \emph{three} and is \emph{cooriented}. The fact that it is an algebraic set makes it a cycle, the fact that it is an hypersurface whose singularities have codimension three allows to define intersection number with it and the fact that is cooriented makes this intersection number an integer. Here coorientation means that $\Sigma$ is two-sided in $L(n)$, i.e. there is a canonical orientation of its normal bundle along its smooth points. Using the above diffeomorphism $L(n)\simeq U(n)/O(n)$ it is easy to choose a positive normal at a smooth point $\Lambda \in \Sigma$: we represent $\Lambda$ as $[M]$ for a unitary matrix $M$ and we take the velocity vector in zero of the curve $t\mapsto [e^{it} M].$

\bex[The train in $L(2)$]
We have seen that $L(2)$ is diffeomorphic to a quadric of signature $(2,3)$ in $\RP^4$; thus it is double covered by $S^1\times S^2$ (i.e. the set of points in $\R^5$ satisfying the equation $x_0^2+x_1^2=x_2^2+x_3^2+x_4^2$ and of norm one).

We fix now a plane $\Delta$ and study the geometry of the train $\Sigma_\Delta$.  We let $\Pi$ be a Lagrangian complement to $\Delta$ and using symmetric matrices chart on $\Pi^{\pitchfork}$ we have:
$$\Sigma_\Delta\cap \Pi^{\pitchfork}\simeq \{S\,|\, \det(S)=0\}.$$
The set of symmetric matrices with determinant zero is a quadratic cone in $\R^3$ with singular point at the origin; to get $\Sigma_\Delta$ we have to add its limit points in $L(2)$ and this results into an identification of the two boundaries components of such a cone. What we get is a Klein bottle with one cycle collapsed to a point.
\eex
The main idea of this paper is to study \emph{generic} maps $f:X\to L(n)$, for $X$ a smooth manifold, and the geometry of the preimage under $f$ of the cycle $\Sigma$ (together with its coorientation). Such a preimage $f^{-1}(\Sigma)$ is what we will call the \emph{induced} Maslov cycle. Sometimes in the sequel the map $f$ will be defined only locally but it will still produce a \emph{Maslov type} cycle on $X$. Our goal is to describe a kind of \emph{duality} that allows to get geometric information on the map $f$  by replacing its study with the one of the geometry of $f^{-1}(\Sigma)$. We will discuss these ideas in greater detail in the next section.  

\bex[Generic loops] Consider a smooth map:
$$\gamma:S^1 \to L(n)$$
transversal to the smooth points of $\Sigma.$ Such a property is generic and we might ask for the meaning of the number of points in $\gamma^{-1}(\Sigma)$. Since the intersection number with $\Sigma$ computes the integer $[\gamma]\in \pi_1(L(n))$, in a very rough way we can write:
\begin{equation}\label{fundgroup}|[\gamma]|\leq b(\gamma^{-1}(\Sigma)),\end{equation}
where the r.h.s. denotes the sum of the Betti numbers, which in this case coincides with the number of connected components (i.e. number of points). This inequality is simply what we obtain by  forgetting the coorientation in the sum defining the intersection number. The comparison through the inequality between what appears on the l.h.s. and what on the r.h.s. is the first mirror of the mentioned duality between the geometric properties of $\gamma$ and the topological ones of $\gamma^{-1}(\Sigma)$.

\begin{remark}[Schubert varieties]\label{remark:schubert}
It is indeed possible to give $L(n)$ a cellular structure using \emph{Schubert varieties} in a fashion similar to the ordinary Grassmannian: the cells are in one to one correspondence with \emph{symmetric} Young diagrams; given one of such diagram the corresponding Schubert cell is the one obtained by considering a flag that is \emph{isotropic} with respect to the symplectic form. More precisely let $\{0\}\subset V_1\subset V_2\cdots\subset V_{2n}=\R^{2n}$ be a complete flag such that $\sigma( V_j, V_{2n-j})=0$ for every $j=1, \ldots, n$ (this means the flag is isotropic; in particular $V_{n}$ is Lagrangian). If now we let $a$ be the partition $a: n\geq a_1\geq a_2\geq \cdots \geq a_n\geq 0$, then the corresponding Schubert variety is:
$$Y_a=\{\Lambda \in L(n)\,|\, \dim (\Lambda \cap V_{n+i-a_i})\geq i\quad \textrm{for}\quad  i=1, \ldots, n\}.$$ The codimension of $Y_a$ is $(|a|+l(a))/2$, where $l(a)$ is the number of boxes on the main diagonal of the associated Young diagram (such a diagram has $a_i$ boxes in its $i$-th row). Since this diagram must be symmetric along its diagonal we see that there are only $2^n$ possible good partitions (see \cite{HeinSottileZelenko} for more details on this approach). Geometrically this shows that the combinatorics of the cell structure of the Grassmannian $G(n,2n)$ descends (by intersection) to the one of $L(n).$ Moreover, since the incidence maps have even degree, cellular homology with $\Z_2$ coefficients gives again the above formula for $H^*(L(n)).$

Notice in particular that $\Sigma$ is a Schubert variety: letting the $n$-th element of the isotropic flag to be $\Delta$ itself, then:
$$\Sigma_\Delta=\{\Lambda \in L(n)\,|\, \dim (\Lambda \cap \Delta)\geq 1\}=Y_{(1, 0, \ldots, 0)}.$$

\end{remark}
\begin{example}[Schubert varieties of $L(2)$]
We consider again the case of $L(2)$ and fix an isotropic flag $\{0\}\subset V_1\subset \Delta \subset V_3\subset \R^4.$ 
The cell structure is given by the four following possible partitions $(0,0), (1,0) ,(2,1) ,(2,2).$ Let us see how the corresponding Schubert varieties look like. To this end let us write $\R^{2n}=\Delta \oplus \Pi$, where $\Pi$ is a Lagrangian complement to $\Delta.$ In this way every $\Lambda$ in $\Pi^{\pitchfork}$ is of the form $\Lambda=\{(x, Sx)\,|\, S=S^{T}\}$.

We immediately get $Y_{(0,0)}=L(n);$ moreover we have already seen that $\Sigma_\Delta=Y_{(1,0)}.$ The Schubert variety $Y_{(2,2)}$ equals $\Delta$ itself (in the symmetric matrices coordinates it is the zero matrix).

Finally we have $Y_{(2,1)}=\{\Lambda \,|\, \Lambda \supset V_1,\, \Lambda \subset V_3\}.$ The intersection of this variety with $\Pi^{\pitchfork}$ equals all the symmetric matrices $S$ whose kernel contains $V_1\subset \Delta$: such matrices are all multiple one of the other and they form a line, thus $Y_{(2,1)}\simeq \RP^1$.\\

\end{example}
\eex \section{Pencils of real quadrics}

\subsection{Local geometry and induced Maslov cylces}
In this section we study in more detail the local geometry of the Lagrangian Grassmannian. If no data are specified, being a differentiable manifold, each one of its points looks exactly like the others. Once we fix one of them, say $\Delta$, the situation drastically enriches: we have seen, for example, that we can choose a cycle $\Sigma_\Delta$ representing the generator of the first cohomology group.

The following proposition gives a more precise structure of the local geometry we obtain on $L(n)$ after we have fixed one of its points $\Delta.$
\bp\label{stratification}Let $\Delta$ in $L(n)$ be fixed. Every $\Lambda\in L(n)$ has a neighborood $U$ and a smooth algebraic \emph{submersion}:
$$\phi:U \to W,$$
where $W$ is an open set of the space of quadratic forms on $ \Delta\cap \Lambda\simeq \R^k$, satisfying the following properties:
\bi
\iii[1.] $(d_{\Lambda} \phi )\dot{\Lambda}=\dot{\Lambda}|_{\Delta \cap \Lambda};$
\iii[2.] $\dim(\ker \phi(\Pi))=k-\dim(\Delta \cap \Pi)$ for every $\Pi$ in $U.$
\iii[3.] for every $\Pi$ in $W$ the fiber $\phi^{-1}(\Pi)$ is contractible.\ei
\ep

Let $\Delta'$ be a lagrangian complement to $\Delta$ transversal to $\Lambda.$ Then, giving coordinates to the open set $\{\Pi\in L(n)\,|\,\Pi\pitchfork \Delta'\}$ using symmetric matrices, the Proposition is just a reformulation of Lemma 2 from \cite{Agrachev1}. 

The fact that $\phi$ is a \emph{submersion} allows to reduce the study of properties of $L(n)$ to smaller Grassmannians, via the Implicit Function Theorem. For the first property, recall that we have a natural identification of the vector space $T_\Lambda L(n)$ with the space of quadratic forms on $\Lambda$; each one of these quadratic forms can be restricted to the subspace $\Delta \cap \Lambda$ and this restriction operation is what $d_\Lambda\phi$ does. The second property says that $\phi$ transforms the combinatorics of intersections with $\Delta$ with the one of the kernels of the corresponding quadratic forms.

Thus locally $\Sigma_\Delta$ looks like the space of degenerate quadratic forms and it is interesting to see how all these local charts are glued together. Let us consider a $\Lambda$ in $\Sigma_\Delta$ and some $\Pi_1$ Lagrangian complement to $\Delta$ such that $\Pi_1\pitchfork\Lambda.$ Given a symplectyc transformation $\psi:\R^{2n}\to \R^{2n}$ preserving $\Delta$,
the matrix $T$ representing it in the coordinates given by the Lagrangian splitting $\Delta\oplus \Pi_1$ has the form: 
$$T=\left(\begin{matrix}A^{-1}&BA^T\\0&A^T\end{matrix}\right)\quad \textrm{with}\quad B=B^T.$$
If $\Lambda$ is represented by the symmetric matrix $S$, the change of coordinates $\psi$ changes its representative to $(A^TSA)(I+BA^TSA)^{-1}$ (indeed this formula works for every $\Lambda$ transversal to $\Pi$).

\begin{remark}[Local topology of the train]
The local topology of $\Sigma_\Delta$ can be described using Proposition \ref{stratification}. Let $B_\Lambda$ be a small ball centered at $\Lambda \in \Sigma_\Delta$ with boundary $S_\Lambda =\partial B_\Lambda$. Then the intersection $B_\Lambda \cap \Sigma_\Delta$ is contractible: it is a cone over the intersection $S_\Lambda \cap \Sigma_\Delta;$ moreover $S_\Lambda \cap\Sigma_\Delta$ is Spanier-Whitehead dual to a union of ordinary Grassmannians and:
\begin{equation}\label{discriminant}H^*(S_\Lambda \cap \Sigma_\Delta)\simeq \bigoplus_{j=0}^{k} H_*(G(j,k)),\quad k=\dim (\Lambda \cap \Delta)\end{equation}
Theorem 3 from \cite{complexity} gives the statement for $\Lambda=\Delta$ and the general result follows by applying Proposition \ref{stratification}.
\end{remark}

For every $r\geq 1$ we can define the sets:
$$\Sigma_\Delta^{(r)}=\{\Lambda \in L(n)\,|\, \dim (\Lambda \cap \Delta)\geq r\}\quad\textrm{and}\quad Z_r=\Sigma_\Delta^{(r)}\backslash \Sigma_\Delta^{(r+1)}.$$
Using this notation, Proposition \ref{stratification} implies that  $\Sigma_\Delta$ is Whitney stratified by $\bigcup_r Z_r$ and the codimension of each $Z_r$ in $L(n)$ is ${r+1\choose 2}$
(the reader is referred to \cite{BCR} for properties of such stratifications). 

 \begin{remark}[Cooorientation revised]Let $\Lambda $ be a smooth point of $\Sigma_\Delta$ and $\gamma:(-\epsilon, \epsilon)\to L(n)$ be a curve transversal to all strata of $\Sigma_\Delta$ and with $\gamma(0)=\Lambda$ (the transversality condition ensures that $\gamma$ meets only $\Sigma_\Delta \backslash \Sigma_\Delta^{(2)}$, i.e. the set of smooth points of $\Sigma$). Since $T_\Lambda L(n)\simeq Q(\Lambda)$,  the velocity $\dot{\gamma}(0)$ can be interpreted (by restriction) as a quadratic form on $\Lambda \cap \Delta$. Proposition \ref{stratification} together with the transversality condition ensures that this restriction is nonzero. We say that the curve $\gamma$ is positively oriented at zero if $\dot{\gamma}(0)|_{\Lambda \cap \Delta}>0.$ Since this definition is intrinsic, it gives a coorientation on $\Sigma$ and it is not difficult to show that it coincides with the above given one.\end{remark}

\bdeff[Induced Maslov cycle] Let $X$ be a smooth manifold and $f:X\to L(n)$ be a map transversal to all strata of $\Sigma=\Sigma_\Delta$. The cooriented preimage $f^{-1}(\Sigma)$ will be called the Maslov cycle induced by $f$.\edeff

A generic map $f:X\to L(n)$ is indeed transversal to all strata of $\Sigma$ and $f^{-1}(\Sigma)$ is itself Nash stratified (its strata being the preimage of the strata of $\Sigma$); the transversality condition ensures that the the normal bundle of the smooth points of $f^{-1}(\Sigma)$ (which is the pull-back of the normal bundle of $\Sigma$) has a nonvanishing section, i.e. the induced Maslov cycle also has a coorientation.

\subsection{Pencils of quadrics}
We turn now to the above mentioned duality between the geometry of a map $f:X\to L(n)$ transversal to all strata of $\Sigma_\Delta$ and the cooriented cycle induced by $f.$ We consider a specific example, namely the case of a map from the sphere, whose image is contained in one coordinate chart.\\
More precisely let $\Delta \oplus \Pi\simeq \R^{2n}$ be a Lagrangian splitting and $W\simeq \R^k$ be a linear subspace of $\Pi^{\pitchfork}\simeq Q(\Delta)$ (the space of quadratic forms on $\Delta$):
$$W=\textrm{span}\{q_1, \ldots, q_k\}\quad\textrm{with}\quad q_1,\ldots,q_k\in Q(\Lambda)\simeq Q(n)$$
(here $Q(n)$ denotes the space of quadratic forms in $n$ variables).\\
Notice that the above isomorphism is defined once a scalar product on $\Delta$ is given: this allows to identify symmetric matrices with quadratic forms.

In this context $W$ is called a \emph{pencil} of real quadrics; the inclusion $S^{k-1}\hookrightarrow W$ defines a map: 
$$f:S^{k-1}\to Q(n)$$
 and for a generic choice of $W$ such a map is transversal to all strata of $\Sigma=\Sigma_\Delta$. Notice that $\Sigma$ equals the discriminant of the set of quadratic forms in $n$ variables and equation (\ref{discriminant}) gives a descritpion of its cohomology.

To every linear space $W$ as above we can associate an algebraic subset $X_W$ of the real projective space $\RP^{n-1}=\mathbb{P}(\Delta)$ (usually referred to by algebraic geometers as the \emph{base locus} of $W$): 
$$X_W=\{[x]\in \RP^{n-1}\, |\, q_1(x)=\cdots=q_k(x)=0\}.$$
The study of the topology of $X_W$ was started by Agrachev  in \cite{Agrachev1, Agrachev2} and continued by Agrachev and the second author in \cite{AgrachevLerario}. 

\begin{remark}[The spectral sequence approach] The main idea of Agrachev's approach is to study the Lebesgue sets of the \emph{positive inertia index} function on $W,$ i.e. the number of positive eigenvalues $\textrm{i}^+(q)$ of a symmetric matrix representing $q$.  More specifically we can consider:
$$W^{j}=\{q\in W\, |\, \textrm{i}^+(q)\geq j\},\quad j\geq 1,$$
and Theorem A from \cite{AgrachevLerario} says that roughly we can take the homology of these sets as the homology of $X_W$:
$$\bigoplus_{j=1}^nH^*(W, W^j)\quad \textrm{``approximates"} \quad H^*(X_W).$$ 
The cohomology classes from $H^*(W,W^j)$ are just the \emph{canditates} for the homology of $X_W$. The requirements they have to fulfill in order to represent effective classes in $H^*(X_W)$ are algebro-topological conditions. The way to make these statements precise is to use the language of spectral sequences (the above conditions on the canditates translate into them being in the kernels of the {differentials} of the spectral sequence). The reader is referred to \cite{AgrachevLerario} for a detailed treatment.  \end{remark}
Going back to the map $f:S^k\to Q(n)$ defined by $W,$ for simplicity of notation we will set:
$$\Sigma_W^{(r)}=S^{k-1}\cap \Sigma^{(r)}.$$
Thus to all these data there correspond two objects: $X_W\subset \RP^{n-1}$ and $\Sigma_W^{(1)}\subset S^{k-1}.$
The induced Maslov cycle is $\Sigma_W$: notice that the cohomology class it represents in $H^1(S^{k-1})$ is clearly zero, though its geometry has a nontrivial meaning.
In fact we can relate the sum of the Betti numbers of $X_W$ to the ones of $\Sigma_W$ and its singular points:
\begin{equation}\label{duality}b(X_W)\leq n+\frac{1}{2}\sum_{r\geq 1}b(\Sigma_W^{(r)})\quad \textrm{for a generic $W$.}\end{equation}
This formula is one of the expressions of the above mentioned duality: the l.h.s. is the \emph{homological complexity} of the intersection of $k$ quadrics in $\RP^{n-1}$, the r.h.s. is the complexity of the Maslov cycle induced on the span of these quadrics. The reader should compare (\ref{duality}) with (\ref{fundgroup}): in both cases the complexity of the induced Maslov cycle gives a restriction (in the form of an upper bound) on some geometrical invariants associated to $f$.

\begin{example}[The intersection of three quadrics]
Let us consider the intersection $X$ of \emph{three} quadrics in $\RP^{n-1}$. Such intersection arises by considering a three dimensional space $W=\textrm{span}\{Q_1, Q_2, Q_3\}$ in a coordinate chart $\Pi^\pitchfork\simeq Q(\Delta).$ Hence the $Q_i$ are symmetric matrices and $X$ is given by the equations $q_1=q_2=q_3=0$ on $\mathbb{P}(\Delta)$; notice that the definition of each $q_i$ depends on the choice of a scalar product on $\Delta$ but two different choices give the same $X$ up to a projective equivalence. The induced Maslov cycle is the curve $\Sigma_W$ on $S^2$ given by the equation:
$$\det( x_1 Q_1+x_2 Q_2+x_3 Q_3)=0,\quad (x_1,x_2,x_3)\in S^2\subset W.$$
This is a degree $n$ curve on $S^2$ and for a generic choice of $W$ it is smooth: in fact $\Sigma_W=S^2\cap \Sigma_\Delta$ and since the codimension of $\textrm{Sing}(\Sigma_\Delta)$ is \emph{three}, by slightly perturbing $W$ this singular locus can be avoided on the sphere.\\
The curve $\Sigma_W$ has at most $O(n^2)$ components and the manifold $X$ at most $O(n^2)$ ``holes" (the sum of its Betti numbers is less than $n^2+O(n)$) ; in this case equation (\ref{duality}) tells that:
$$|b(X)-b_0(\Sigma)|\leq O(n),$$
i.e. if we replace the homology of $X$ with the one of the associated Maslov cycle the error of such replacement has order $O(n).$ The coorientation of the induced Maslov cycle in this case assigns a number $\pm 1$ to each oval of the curve $\Sigma_W$: this number is obtained by looking at the change of the number of positive eigenvalues when crossing the oval. The knowledge of the coorientation on each oval allows to compute the error term in (\ref{duality}); the reader is referred to \cite{Agrachev1, AgrachevLerario, complexity}.
\end{example}

\section{Geometry of Gauss maps}
\subsection{Lagrange submanifolds of $\mathbb{R}^{2n}$}
Consider a Lagrangian submanifold $M$ of the symplectic space $T^*\R^n\simeq \R^{2n}$. The \emph{Gauss map} of $M$ is:
$$\nu:M\to L(n)$$
and associates to each point $x\in M$ the tangent space $T_xM$ (which is by definition a Lagrangian subspace of $\R^{2n}$). 

We consider the Lagrangian splitting $\R^{2n}=\Pi\oplus \Delta$ and we are interested in the description of the induced Maslov cycle $\nu^{-1}(\Sigma_\Delta)$ on $M.$ To this end we consider the projection on the first factor $\pi:\R^{2n}\to \Pi$ and its restriction to $M$:
$$\pi|_M:M\to \Pi.$$
The critical points of $\pi|_M$ are those points $x$ in $M$ such that the tangent space $T_xM$ does not intersect $\Delta$ transversally; in other words:
\begin{equation}\label{eq:critical}\textrm{Crit}(\pi|_M)=\nu^{-1}(\Sigma_\Delta).\end{equation}
Thus the induced Maslov cycle in this case coincides with the set of critical points of a map from $M$ to $\R^n$: this critical set represents the Poincaré dual of $w_1(TM)$, the first Stiefel-Whitney class of $TM$ (see Remark \ref{remark:cc} below). In fact $\nu$ pulls-back the tautological bundle $\tau(n)$ of $L(n)$ to the tangent bundle of $M$ and, by functoriality of characteristic classes, it also pulls-back the first Stiefel-Whitney class of $\tau(n)$ to $w_1(TM).$ Notice that $w_1(\tau(n))$ equals the modulo two reduction of the universal Maslov class $\mu$ defined above.

\begin{remark}[Characteristic classes revised]\label{remark:cc} Consider an $n$-dimensional manifold $M$ and a smooth function $f:M\to\R^{n-k+1}.$
For a generic $f$ we can relate the $k$-th Stiefel-Whitney class of $M$ to the critical points of $f$ by:
\begin{equation}\label{cc}w_{k}(TM)=\textrm{Poincaré dual of } \textrm{Crit}(f).\end{equation}

For $k=n$ the generic $f$ is a Morse function and  $w_n(TM)\in H^n(M)\simeq \Z_2$ is the Euler characteristic of $M$ modulo two, thus the previous equations reads $\chi(M)\equiv \textrm{Card}(\textrm{Crit}(f)) \textrm{mod }2.$

In the case $k=1$ we can apply (\ref{cc}) to: 
$$f=\pi|_M:M\to \R^n,$$
 and equation (\ref{eq:critical}) implies that the Maslov cycle induced by $\nu$ represents the Poincaré dual of $w_1(TM).$
 
We know from Remark \ref{remark:schubert} that the cohomology of $L(n)$ is generated by the Poincaré duals of its Schubert varieties. Each of these varieties is labelled using \emph{symmetric Young diagrams} and their intersections are computed using Schubert calculus. The variety corresponding to the diagram having only one box is $Y_{(1, 0, \ldots, 0)}$: this is the train of $\Delta$ (the middle space in the isotropic flag) and it represents the Poincaré dual of $\mu=w_1(\tau(n))
$ (again reduction modulo two is considered).

\end{remark}
\begin{example}[Surfaces in $\R^4$] Among compact orientable surfaces $S$, the only one that admits a Lagrangian embedding into $\R^4$ is the torus (in particular $w_1(TS)=\nu^*w_1(\tau(n))=0$). If $S$ is nonorientable and is not the Klein bottle, it has a Lagrangian embedding in $\R^4$ if and only if:
$$\chi(S)\equiv 0 \mod 4.$$
The Klein bottle cannot be realized as a Lagrangian submanifold of $\R^4.$\end{example}

\begin{example}[Spheres] No exact Lagrangian embeddings (i.e. embeddings for which the two-form $\sigma$ pulls-back to an exact one-form) exist in $\R^{2n}.$ In particular simply connected manifolds, as the spheres, cannot be Lagrangian submanifolds of $\R^{2n}.$ \\
If we only require the differential to be injective, a Lagrangian \emph{immersion} $j: S^n\to \R^{2n}\simeq \C^n$ is given by:
$$j:(x,y)\mapsto (1+2iy)x,\quad \textrm{where}\quad S^n=\{(x,y)\in \R^n\times \R\,|\, \|x\|^2+y^2=1\}.$$
In the case $n=1$ the image of $j$ is an eight-shaped curve; this immersion fails to be injective at the north and the south pole only.

\end{example}

\subsection{Lagrangian maps}
Generalizing the construction of the previous section, we consider $M$ Lagrangian submanifold  of the symplectic manifold $T^*N$ (with the standard symplectic structure); we denote by $\pi:T^*N\to N$ the bundle projection. In this case we do not have a global Gauss map, but in analogy with (\ref{eq:critical}) we can still define the induced Maslov cycle as:
$$\Sigma_M=\textrm{Crit}(\pi|_M).$$
The case of a submanifold $M$ of $T^*N$ projecting to $N$ is itself a special case of a \emph{Lagrangian map}; this is defined as follows. First we say that a fibration $\pi:E\to N$ is Lagrangian if $E$ is a symplectic manifold and each fiber is Lagrangian. A Lagrangian map is thus a smooth map $f:M\to N$  between manifolds of the same dimension obtained by composition of a Lagrangian inclusion $i:M\to E$ followed by $\pi:$
$$f:M\stackrel{i}{\longrightarrow}E\stackrel{\pi}{\longrightarrow}N.$$ 
We refer the reader to \cite{ArnoldMaMeClMe} for more details and examples.

\begin{example}[Normal Gauss maps of hypersurfaces]
Consider a smooth oriented hypersurface $M$ in $\R^{n+1}$; the \emph{normal Gauss map} of $M$ is the map:
$$f:M\to S^n,\quad x\mapsto \textrm{oriented normal of $M$ at $x$}.$$
This map is Lagrangian; in fact we can set $E=T^*S^n\simeq TS^n$ with projection $\pi:E\to S^n$ and define the Lagrangian inclusion $i:M\to E$ as $x\mapsto (f(x), \textrm{proj}_{T_xM}x).$
The image in $S^n$ of the induced Maslov cycle under $f$ is called the \emph{focal surface} of $M.$

\end{example}

Thus a Lagrangian map $f:M\to N$ is a special case of map between two manifolds of the same dimension. Proposition \ref{stratification} allows to give a local description of the set of critical points of a Lagrangian map: it is a \emph{cooriented} hypersurface in $M$, smooth outside a set of codimension three. The set of critical values of a Lagrangian map $f$ is called a \emph{caustic}.

\section{Lagrange multipliers}
Let $\U$ be an open set in a Hilbert space (or a finite dimensional manifold) and let $M$ be a smooth $n$-dimensional manifold. Assume we have a pair of  smooth maps 
$\FF: \U\to M,$ and $\JJ:\U\to \R.$ We want to characterize critical points of the functional $\JJ$ when restricted to level sets of $\FF$:
\bqn \label{eq:crit}
\min_{\FF^{-1}(x)} \JJ, \qquad x\in M.
\eqn
\index{critical poinr!constrained}
Recall that for a smooth function $f:M \to \R$ and a smooth submanifold $N\subset M$ a point $x\in N$ is said a \emph{critical point} of $f\big|_{N}$ if $d_{x}f\big|_{T_{x}N}=0$.
We state the geometric version of the Lagrange multipliers rule, which characterizes regular constrained critical points.
\bp[Lagrange multipliers rule] \label{p:lmr} Assume $u\in \U$ is a regular point of $\FF:\U \to M$ such that $\FF(u)=x$. Then $u$ is a critical point of $\JJ\big|_{\FF^{-1}(x)}$ if and only if:
\bqn \label{eq:lagrm}
\exists \lambda \in T^{*}_{x}M,\  \lambda \neq 0, \qquad \text{s.t. } \qquad d_{u}\JJ=\lambda \,D_{u}\FF.
\eqn
\ep
The above discussion suggests to consider the pairs $(u,\lambda)$ such that the identity $d_{u}\JJ=\lambda \,D_{u}\FF$ holds true. More precisely we should consider the pair $(u,\lambda)$ as an element of the pullback bundle $\FF^{*}(T^{*}M)$, and set
$$\Crit_{\FF,\JJ}=\{(u,\lambda)\in \FF^{*}(T^{*}M) | \, d_{u}\JJ=\lambda \,D_{u}\FF\}$$
Notice that by definition of pullback bundle, if $(u,\lambda)\in \FF^{*}(T^{*}M)$, then $\FF(u)=\pi(\lambda)$ ($\pi:T^*M\to M$ is the bundle projection).
The study of the geometry of the set $\Crit_{\FF,\JJ}$ leads us to investigate the constrained critical points for the whole family of problems \eqref{eq:crit}, as $x$ varies on $M$.
The following regularity condition ensures that $\Crit_{\FF,\JJ}$ has nice properties: the pair $(\FF,\JJ)$ is said to be a \emph{Morse pair} (or a \emph{Morse problem}) if the function 
\bqn\label{eq:theta}
\theta: \FF^{*}(T^{*}M) \to T^{*}\U, \qquad (u,\lambda) \mapsto d_{u}\JJ- \lambda\, D_{u}\FF.
\eqn
is transversal to the zero section in $T^{*}\U$.
Notice that, if $M=\{0\}$, then $\FF$ is the trivial map and with this definition we have that $(\FF,\JJ)$ is a Morse pair if and only if $\JJ$ is a Morse function. 

If $(\FF,\JJ)$ defines a Morse problem, then $\Crit_{\FF,\JJ}$ is a smooth $n$-dimensional manifold in  $\FF^{*}(T^{*}M)$. In the case when $\U$ is a finite dimensional manifold this is easy to show it, since by a standard transversality argument:
\begin{align*}
\dim\, \Crit_{\FF,\JJ}&=
\dim\, \FF^*(T^*M)-\dim\,\U\\
&=(\dim\, \U+ \tx{rank}\, T^*M) -\dim\, \U\\
&= \tx{rank}\, T^*M =n
\end{align*} 
The above argument is no more valid in the infinite dimensional case but one can show that the same result holds (under some additional technical assumptions).

Let us now consider the map $\overline{\FF}: F^{*}(T^{*}M)\to T^{*}M$ given by $(u,\lambda)\mapsto \lambda$. We can consider the set $\overline{C}_{F,J}=\overline{\FF}(\Crit_{\FF,\JJ})$ in $T^*M:$
\bqn\label{eq:diagr10}
\xymatrix{
\Crit_{\FF,\JJ} \ar[d]_{\overline{\pi}} \ar[r]^{\overline{\FF}}
& T^{*}M \ar[d]^{\pi} \\
\U \ar[r]_{\FF}& M }
\eqn

It turns out that $\overline{\FF}$ is an exact Lagrangian immersion, i.e. it pulls-back the Liouville form $p\,dq$ to an exact form. 

We assume now that $\Lagr_{\FF,\JJ}$ is an embedded submanifold (and not only immersed). 
\bp \label{p:badpr}
Let $(\FF,\JJ)$ be a Morse problem and assume $(u,\lambda)$ is a Lagrange multiplier such that $u$ is a regular point for $\FF$, where $\FF(u)=x$. The following properties are equivalent:
\bi
\iii[1.] $\emph{\textrm{Hess}}_{u}\, \JJ \big|_{\FF^{-1}(x)}$ is degenerate,
\iii[2.] $(u,\lambda)$ is a critical point for the map $\pi \circ \overline{\FF}: \Crit_{\FF,\JJ} \to M.$
\ei

\ep
We will discuss the proof in a special case in the next section; the general proof follows the same line.

Notice that $\pi \circ \overline{F}:C_{F,J}\to M$ is a Lagrangian map; the induced Maslov cycle $\Sigma_{F,J}$, i.e. the set of critical points of this map, coincides with the set of those $(\lambda, u)$ such that the Hessian of $J|_{F^{-1}(F(u))}$ is degenerate at $u$.

\subsection{Morse functions}
Let us consider two Morse functions $f_{0},f_{1}:M\to \R$ and an homotopy of maps $f_{t}:M\to \R$. Then we define $\U=[0,1]\times M$ and: 
$$F:[0,1]\times M\to \R,\qquad F(t,x)=t$$
$$J:[0,1]\times M\to \R,\qquad J(t,x)=f_{t}(x).$$
We have that
$\JJ \big|_{\FF^{-1}(t)}=f_{t}$ and we can study the critical points of the family of maps $\{f_{t}\}_{t\in[0,1]}$ with the Lagrange multipliers technique.
 If $u=(t,x)$, writing $D_{u}J=(\partial_{t}J, \partial_{x}J)$ and $D_{u}F=(1, 0)$  the Lagrange multipliers rule reads
\begin{equation}\label{eq}\begin{cases}
\lambda= \partial_{t}J(t,x)\\
\partial_{x}J(t,x)=0
\end{cases}
\end{equation}
Namely $\Crit_{\FF,\JJ}$ is the set of $(\lambda,t,x)$ such that \eqref{eq} holds true (
the second identity is equivalent to the fact that $x$ is a critical point of $f_{t}$).
This is a system of $n+1$ equations in a $n+2$-dimensional space and 
$\Crit_{\FF,\JJ}$ defines a 1-dimensional manifold if the problem is Morse, i.e. the linearized system in the variables $(\lambda',t',x')$
\begin{equation}\label{eq22}\begin{cases}
\lambda'= \partial^{2}_{tt}J(t,x) t' +\partial^{2}_{xt}J(t,x) x' \\
\partial^{2}_{tx}J(t,x) t' +\partial^{2}_{xx}J(t,x) x'=0
\end{cases}
\end{equation}
is regular, that means 
$\text{rank}\{\partial^{2}_{tx} f,\partial^{2}_{xx}f  \}=n.$ In particular this condition is satisfied if the function $f_{t}$ is Morse for every $t\in [0,1]$. The tangent space to $\Crit_{\FF,\JJ}$ at the point $(\lambda,t,x)$ is the set of $(\lambda',t',x')$ such that \eqref{eq22} are satisfied. 
\bqn\label{eq:diagr0}
\xymatrix{
(\lambda,t,x) \ar[d]_{\overline{\pi}} \ar[r]^{\overline{\FF}}
& (\lambda,t) \ar[d]^{\pi} \\
(t,x) \ar[r]_{\FF}& t }
\eqn
Hence the point $(\lambda,t,x)\in \Crit_{\FF,\JJ}$ is critical for the map if and only if there exists a nonzero element $(\lambda',t',x')$ such that $\pi \circ \overline{\FF}(\lambda',t',x')=t'=0.$ From \eqref{eq22} it is easy to see that this is equivalent to $x'\neq 0$ and  $\partial^{2}_{xx}J(t,x) x'=0.$

Let now $f_t$ be a generic homotopy between two Morse functions $f_0$ and $f_1$. Then the corresponding pair $(F,J)$ defines a Morse problem and the above discussion holds. Moreover the genericity assumption on the homotopy ensures that if  $f_{t_0}$ has a critical point at $x_0$, the Hessian of $f_{t_0}$ at $x_0$ has a one-dimensional kernel. It is indeed possible to show that  near the point $( t_0, x_0)$ the family $f_t$ can be written in coordinates as:
$$f_t(x)=c_0+x_1^3\pm t x_1\pm x_2^2+\cdots\pm x_n^2, \quad t\in [t_0-\epsilon, t_0+\epsilon].$$
As $t$ passes from $t_0-\epsilon$ to $t_0+\epsilon$ two critical points merge or vanish, according to the sign of $\pm tx$ (see \cite{Milnor}).

The induced Maslov cycle $\Sigma_{F,J}$ in this case consists of those points $(\lambda, t,x)$ on $C_{F,J}$ such that $f_t$ is not a Morse function. If $(\lambda, t, x)$ is in $C_{F,J}\backslash \Sigma_{F,J}$, then in a neighborood $[a,b]$ of $t$ the function $t$ is a coordinate for $C_{F,J}$ and we can  ``follow" the critical point $x(t)$. Property 1 of Proposition \ref{p:badpr} implies that as long as $t$ varies on $[a,b]$, the index of such critical point never changes. The genericity assumption on the homotopy implies that if two critical points merge, their indices \emph{must} differ by one. If $(\lambda(s), t(s), x(s))$ is a parametrization of $C_{F,J}$ near a point $(\lambda(0), t(0), x(0))\in \Sigma_{F,J}$, the change in the \emph{sign} of the determinant of the Hessian of $f_{t(s)}$ at $x(s)$ when passing through $s=0$ is determined by the coorientation of $\Sigma_{F,J}$ at $(\lambda(0), t(0), x(0)).$

In this case the number of points of  $\Sigma_{F,J}$ tells how many functions in our family are not Morse; the coorientation tells how the Morse index changes when two critical points merge or vanish.

\begin{example}[Depth of Morse functions]
Assume $M$ is a smooth hypersurface in $\R^n$ defined by a polynomial of degree $d$ and $p_0, p_1$ are two Morse functions obtained by restricting to $M$ two polynomials of degree $k\geq d$. Using the above technique it is possible to prove that $p_0$ and $p_1$ can be joined by a homotopy $p_t:M\to \R$ such that:
$$\textrm{Card}\{t\in [0,1]\,|\, p_t \textrm{ is not Morse}\}\leq dk^n (d+nk).$$
In the case $k\leq d$ the bound is $d^{n+2}(n+1).$\end{example}

\subsection{Riemannian and sub-Riemannian geometry}\label{s:rg}
In this section we discuss how the problem of finding geodesics in Riemannian or sub-Riemannian geometry
fits in the above setting. For a comprehensive presentation of Riemannian and sub-Riemannian geometry see for instance \cite{nostrolibro,montgomerybook}.

A sub-Riemannian manifold is a triple $(M,\mathcal{D},g)$ where $M$ is a smooth manifold and $\mathcal{D}$ is a constant rank $k\leq n$ distribution endowed with a scalar product $g$ on it. The case $k=n$, i.e. when $\mathcal{D}=TM$, corresponds to Riemannian geometry. 

A curve on $M$ defined on the interval $[0,1]$ is said horizontal if it is almost everywhere tangent to the distribution. Once fixed a local orthonormal basis of vector fields $f_{1},\ldots,f_{k}$ on $\mathcal{D}$, every horizontal curve is described by the dynamical system: 
\bqn \label{eq:ex}
\dot x(t)=\sum_{i=1}^{k}u_{i}(t)f_{i}(x(t)),\qquad x(0)=x_{0},
\eqn
for some choice of the control $u;$ the length of such horizontal curve is defined by:
$$\ell(u)=\int_{0}^{1} \sqrt{g(\dot{x}(t),\dot{x}(t))}\,dt=\int_{0}^{1}\sqrt{\sum_{i=1}^{k}u_{i}^{2}(t)}\,dt.$$
It is well known that the problem of minimizing the length  with fixed final time is equivalent, by Cauchy-Schwartz inequality, to the minimization of the energy 
$$J(u)=\int_{0}^{1}\sum_{i=1}^{k}u_{i}^{2}(t)\,dt.$$
For this reason it is convenient to parametrize horizontal curves by admissible controls $u\in L^{2}([0,1],\R^{k})$. By the classical theory of ODE, for every such control $u$ and every initial condition $x_{0}\in M$, there exists a unique solution $x_u$  to the Cauchy problem
\eqref{eq:ex}, defined for small time (see for instance \cite{agrachevbook} for a proof).

The resulting local flow defined on $M$ by the ODE associated with $u$, i.e. the family of diffeomorphisms $P_{0,t}:M\to M$, defined by $P_{0,t}(x):=x_{u}(t)$ is smooth in the space variable and Lipschitz in the time variable. 
Analogously one can define the flow $P_{s,t}:M\to M$ fixing the initial condition at time $s$, i.e. $x(s)=x_{0}$ ($P_{s,t}$ is defined for $s,t$ close enough).

Fix a point $x_{0}\in M$.  The \emph{end-point map} of the system \eqref{eq:ex} is the map 
$$F: \U\to M,\qquad u\mapsto x_{u}(1),$$
where $\U\subset L^{2}([0,1],\cb)$ is the open subset  of controls $u$ such that the solution $t\mapsto x_{u}(t)$ to the Cauchy problem 
 \eqref{eq:ex}
  exists and is defined on the whole interval $[0,1]$.
  The end-point map is a smooth map. Moreover its differential $D_{u}F: L^{2}([0,1],\R^{k})\to T_{x}M$ at a point $u\in U$ is computed by the following well-known  formula (see \cite{agrachevbook})
\bqn\label{eq:duev}
D_{u}F(v)=\sum_{i=1}^{k}\int_{0}^{1}v_{i}(s)(P_{s,1})_{*}f_{i}(x_{u}(s)) ds,\qquad v\in L^{2}([0,1],\R^{k}).
\eqn
where $x_{u}(t)$ is the trajectory associated with $u$ and $x=x_{u}(1)$. \\
Notice that when $u=0$ we have rank$\,D_{0}F=\text{rank} \,\mathcal{D}=k$. Indeed $x_{u}(t)\equiv x_{0}$ and the above formula reduces to: 
$$D_{0}F(v)=\sum_{i=1}^{k}\alpha_{i}f_{i}(x_{0}), \qquad \alpha_{i}=\int_{0}^{1}v_{i}(s) \,ds.$$


In this framework, the problem of finding constrained critical points of the functional $J:U\to \R$ on the level set $F^{-1}(x)$ is equivalent to find critical points of the energy among those curves that join $x_{0}$ to $x$ in fixed final time equal to 1. 

Hence the solutions of the problem \eqref{eq:crit} represent exactly sub-Riemannian \emph{geodesics} starting at $x_0$ and ending at $x$. 

Notice that in the Riemannian case the map $F$ is always a submersion, while in the sub-Riemannian case it can happen that $\text{rank}(D_{u}F)<n$ for some $u$ (this is the case for the control $u=0$ as we explained above). In this case $u$ is said \emph{abnormal} and $x_{u}$ is an \emph{abnormal geodesic}. If $u$ satisfies the Lagrange multipliers rule $\lambda D_{u}F=D_{u}J$ for some $\lambda$, then $u$ is said \emph{normal} and $x_{u}$ is a \emph{normal geodesic} (this happens in particular at regular point of $F$). A control $u$ can be at the same time normal and abnormal.

In what follows we focus our attention to \emph{strongly normal} controls, i.e. those controls such that all the family  $u_{s}(t):= s u(st)$ is not abnormal for all $s\in]0,1]$. Notice that, by the linearity of \eqref{eq:ex} with respect to $u$, we have $x_{u_{s}}(t)=x_{u}(st)$. Notice also that in Riemannian geometry all geodesics are strongly normal.

Given a sub-Riemannian structure on a manifold $M$ it is natural to build the \emph{sub-Riemannian Hamiltonian} $H:T^{*}M\to \R$ defined by
$$H(\lambda)=\frac{1}{2}\|\lambda\|^{2},\qquad \|\lambda\|=\sup_{v\in \mathcal{D}_{q},|v|\leq 1} |\langle\lambda,v\rangle|.
$$
This is a smooth function on $T^{*}M$ which is quadratic on fibers. The canonical symplectic structure allows to define a vector field $\overrightarrow{H}$ by the identity
$\sigma(\cdot, \overrightarrow{H})=dH.$
The flow of $\overrightarrow{H}$ defines the normal geodesic flow and characterizes the manifold of Lagrange multipliers as follows.
\bp The sub-Riemannian pair $(F,J)$ defines a Morse problem. Moreover the manifold of Lagrange multipliers satisfies $\overline{C}_{F,J}=e^{\overrightarrow{H}}(T^{*}_{x_{0}}M)$.
\ep

We discuss some related ideas, giving an outline of the proof.

 Let $x\in M$ and  $(u, \lambda)\in \Crit_{F,J}$ associated with a critical point of $J\big|_{F^{-1}(x)}$. Then for every $v\in \ker D_{u}F:$
\begin{align}\label{eq:22}
\text{Hess}_{u} J\big|_{F^{-1}(x)}(v)&=\|v\|_{L^{2}}^{2}-\langle \lambda\, ,\underset{0\leq \tau \leq t \leq 1}{\iint} 
[(P_{\tau,1})_{*}f_{v(\tau)},(P_{t,1})_{*}f_{v(t)}] d\tau dt \rangle.
\end{align}
Indeed one can compute that in coordinates $\text{Hess}_{u} J\big|_{F^{-1}(x)}= D^{2}_{u}J-\lambda D^{2}_{u}F$ and that the second differential of the end-point map is expressed as the commutator
$$D^{2}_{u}F(v,v)=\underset{0\leq \tau \leq t \leq 1}{\iint} 
[(P_{\tau,1})_{*}f_{v(\tau)},(P_{t,1})_{*}f_{v(t)}] d\tau dt, $$
where $P_{s,t}$ is the \emph{non autonomous} flow defined by the control $u$ and $f_{v}=\sum_{i=1}^{k} v_{i} f_{i}$.

Let $(u, \lambda)\in \Crit_{F,J}$. The relation $D_{u}J=\lambda D_{u}F$ can be rewritten as follows, using the fact that $J(u)=\|u\|^{2}_{L^{2}}$:
\bqn\label{eq:LL}
u_{i}(t)=\langle\lambda(t), f_{i}(x(t))\rangle, \qquad \lambda(t):=(P_{t,1})^{*}\lambda\in T^{*}_{x(t)}M.
\eqn
Moreover the curve $\lambda(t)\in T^{*}_{x(t)}M$ is a solution of the Hamiltonian system $\dot{\lambda}(t)=\overrightarrow{H}(\lambda(t))$ and $\lambda(1)=\lambda$.
This allows to parametrize geodesics via their initial covector rather than the final one.

We define the \emph{exponential map} starting from $x_{0}$ as the map:
$$\mathcal{E}: T^{*}_{x_{0}}M\to M,\qquad \mathcal{E}(\lambda_{0})=\pi\circ e^{\overrightarrow{H}}(\lambda_{0}).$$ 
Since $e^{\overrightarrow{H}}(T^*_{x_0}(M))=\overline{C}_{F,J},$ then this map is Lagrangian. Moreover, by homogeneity of the Hamiltonian, for all $t>0$  we have $\mathcal{E}(t\lambda_{0})=\pi\circ e^{t\overrightarrow{H}}(\lambda_{0})=x_{u}(t),$ 
which permits to recover the whole normal geodesic associated with $\lambda_{0}$ (here $u$ is the control defined by \eqref{eq:LL} and $\lambda(t)=(P_{t,0})^{*}\lambda_{0}$). Thus the exponential map parametrizes normal geodesics starting from a fixed point with covectors attached to the fiber $T^{*}
_{x_{0}}M$. If $\lambda_{0}$ is a critical point of $\mathcal{E}$ then the point $x=x_{u}(1)=\mathcal{E}(\lambda_{0})$ is said to be \emph{conjugate} to $x_{0}$ along the geodesic $x_{u}(t)$.
 
\bp Let $x_{u}(t)$ be a strongly normal geodesic joining $x_{0}$ to $x$. The following are equivalent
\bi
\iii[1.] $\emph{\textrm{Hess}}_{u}\, \JJ \big|_{\FF^{-1}(x)}$ is degenerate;
\iii[2.] $x$ is conjugate to $x_{0}$ along $x_{u}(t)$.
\ei
Moreover the geodesic $x_{u}(t)$ loses its local optimality at its first conjugate point.
\ep

By the homogeneity of the Hamiltonian, to study the local optimality of a piece $x_{u}|_{[0,s]}$ of the fixed trajectory $x_{u}$ it is enough to apply the functional $J$ to the control $u_{s}(t)=su(st)$, whose final point is $x_{u_{s}}(1)=x_{u}(s)$.

Thus we have the following picture: the map $\mathcal{E}:T^*_{x_0}M\to M$ is a Lagrangian map with the property that $\mathcal{E}(\lambda_0)$ is the final point of a geodesic $x$ starting at $x_0$; this geodesic is the one associated to the control $u$ defined by equation (\ref{eq:LL}). 

The induced Maslov cycle $\Sigma_{x_0}$, i.e. the set of critical points of $\mathcal{E}$, coincides with the set of those $\lambda\in T^*_{x_0}M$ such that the Hessian of $J|_{F^{-1}(F(u))}$ at the corresponding geodesic is degenerate. We can indeed consider the all ray $\{s\lambda\}_{s>0}$: the image of such ray is the geodesic associated with $\lambda$.  For small $s>0$ the Hessian $\textrm{Hess}_{u_{s}}\, \JJ \big|_{\FF^{-1}(x_{u}(s))}$ is positive definite (as a consequence of formula \eqref{eq:22}), and it becomes degenerate exactly when $s\lambda$ belongs to the induced Maslov cycle $\Sigma_{x_0}$ (in particular the first degeneracy point coincide with the first conjugate point).
 
With  a normal geodesic $x(t)$ (with lift $\lambda(t)$) one can associate also the so-called \emph{Jacobi curve}: 
$$\Lambda(t)=e^{-t\overrightarrow{H}}_{*}T_{\lambda(t)}(T^{*}_{x(t)}M),$$
which is a curve of Lagrangian subspaces in the symplectic space $T_{\lambda_{0}}(T^{*}_{x_{0}}M)$. Using this curve we can compute the index of the Hessian: in fact if $c$ is a cocycle representing the Maslov class $\mu$, we have:
$$c(\Lambda(s))=-\textrm{Ind}\,\textrm{Hess}_{u_{s}}\, \JJ \big|_{\FF^{-1}(x_{u}(s))}.$$

\begin{remark}[General variational problem] 
The results obtained in Section \ref{s:rg} can be extended to the more general case of a non-linear control problem
$$\dot x(t)=f(x(t),u(t)),\qquad x(0)=x_{0},$$
with a Tonelli type integral cost
$$J(u)=\int_{0}^{1}L(x_{u}(t),u(t))dt.$$
Namely we require that $L(x,\cdot)$ is strictly convex and super linear with respect to $u$. 

Under these assumptions it is still possible to characterize the manifold of Lagrange multipliers via the Hamiltonian associated to this problem (see \cite{cime}). 
\end{remark}

\vspace{0.2cm}
{\bf Acknowledgements.} The first author has been supported by the European Research Council, ERC StG 2009 ``GeCoMethods", contract number 239748, by the ANR Project GCM, program ``Blanche", project number NT09-504490.

\input{referenc}

\end{document}

%% file: referenc.tex
%
%
%

%% file: Z-Maslov-nospringer.bbl
\begin{thebibliography}{99.}%
%
%

\bibitem{Agrachev1} A. A. Agrachev: \emph{Topology of quadratic maps and hessians of smooth maps}, Translated in J. Soviet Math. 49, 1990, no. 3, 990-1013.
      
      \bibitem{cime} A. A. Agrachev, \emph{Geometry of optimal control problems and Hamiltonian systems}. Nonlinear and optimal control theory, 1–59, 
Lecture Notes in Math., 1932, Springer, Berlin, 2008. 


\bibitem{agrachevbook} A. A. Agrachev and Y. Sachkov: {\em Control theory from the geometric
  viewpoint}, vol.~87 of Encyclopaedia of Mathematical Sciences,
  Springer-Verlag, Berlin, 2004.

\bibitem{nostrolibro}
A.~Agrachev, D.~Barilari, and U.~Boscain: {\em Introduction to
  {R}iemannian and sub-{R}iemannian geometry ({L}ecture {N}otes),
  http://people.sissa.it/agrachev/agrachev\_files/notes.html},  (2012).
            
 \bibitem{Agrachev2} A. A. Agrachev and R. V. Gamkrelidze: \emph{Quadratic maps and smooth vector valued functions; Euler Characteristics of level sets}, Itogi nauki.
 \bibitem{AgrachevLerario} A. A. Agrachev, A. Lerario: \emph{Systems of quadratic inequalities}, Proc. London Math. Soc. 2012, 105 (3). 
 
 \bibitem{AgrachevQuHo} A. A. Agrachev: \emph{Quadratic Homology}, preprint.
 
  \bibitem{AgCurv} A. A. Agrachev, D. Barilari, P.W.Y. Lee and L.Rizzi: \emph{Curvature for affine control sysyems and sub-Riemannian geometry}, preprint.
 
 
 \bibitem{ArnoldChCl} V. I. Arnold: \emph{On a characteristic class entering into conditions of quantisation}, English translation. Functional Analysis and Its Applications 1, 1-14 (1967)
 
 \bibitem{ArnoldMaMeClMe} V. I. Arnold: \emph{Mathematical Methods of Classical Mechanics}, Springer-Verlag (1989).
 
 \bibitem{ArnoldGivental} V. I. Arnold, A. B. Givental: \emph{Symplectic geometry}, Translated from 1985 Russian original. in Dynamical Systems IV, Encycl. of Math. Sciences 4, Springer, 1-136 (1990)
    
  \bibitem{montgomerybook}
R.~Montgomery:  {\em A tour of subriemannian geometries, their geodesics
  and applications}, vol.~91 of Mathematical Surveys and Monographs, AMS, Providence, RI, 2002.
  
\bibitem{BCR} J. Bochnak, M. Coste, M-F. Roy: \emph{Real Algebraic Geometry}, Springer-Verlag (1998). 

 \bibitem{Fuks} D. B. Fuchs, O. Ya. Viro: \emph{Classical manifolds}, in Topology II, Encyclopaedia of Mathematical Sciences, volume 24, Springer (2000).

  \bibitem{Hatcher} A. Hatcher: \emph{Algebraic Topology}, Cambridge University Press (2002).
 
\bibitem{HeinSottileZelenko} N. Hein, F. Sottile, I. Zelenko: \emph{A congruence modulo four in real Schubert calculus}, 	arXiv:1211.7160

\bibitem{complexity} A. Lerario, \emph{Complexity of intersection of real quadrics and topology of symmetric determinantal varieties},  arXiv:1211.1444

\bibitem{Milnor} J. Milnor: \emph{Lectures on the h-cobordism theorem}, Notes by L.Siebenmann and J.Sondow, Princeton Math. Notes (1965).
%
%
\end{thebibliography}
